\documentclass[psamsfonts,reqno]{amsart}
\usepackage{amscd}

\theoremstyle{plain}

\newtheorem{lemma}{Lemma}[section]
\newtheorem{theorem}{Theorem}

\newtheorem{corollary}[lemma]{Corollary}
\newtheorem{claim}{Claim}

\newtheorem*{stat}{\name}
\newcommand{\name}{testing}

\theoremstyle{definition}
\newtheorem{definition}[lemma]{Definition}

\newtheorem{problem}{Problem}

\theoremstyle{remark}

\newenvironment{all}[1]{\renewcommand{\name}{#1}\begin{stat}}
                        {\end{stat}}

\newcommand{\qedc}{{\qed}~{\rm Claim~{\theclaim}.}}

\newenvironment{cproof}
{\begin{proof}[Proof of Claim.]}
{\qedc\renewcommand{\qed}{}\end{proof}}

\numberwithin{equation}{section}

\def\vv<#1>{\langle#1\rangle}
\newcommand{\set}[1]{\{\,#1\,\}}
\newcommand{\setm}[2]{\set{#1\mid#2}}
\newcommand{\famm}[2]{\left\langle\,#1\mid#2\,\right\rangle}

\newcommand{\two}{\mathbf{2}}
\newcommand{\clp}{\textup{CLP}}

\newcommand{\xa}{\boldsymbol{a}}

\newcommand{\xq}{\boldsymbol{q}}
\newcommand{\xx}{\boldsymbol{x}}
\newcommand{\xy}{\boldsymbol{y}}

\DeclareMathOperator{\Conc}{Con_c}
\DeclareMathOperator{\Id}{Id}

\newcommand{\jz}{$\set{\vee,0}$}
\newcommand{\jzh}{\jz-ho\-mo\-mor\-phism}

\begin{document}

\title[Lifting problems for congruence lattices]%
{Unsolvable one-dimensional lifting problems for congruence
lattices of lattices}

 \author[J.~T\r uma]{Ji\v r\'\i\ T\r uma}
 \address{Department of Algebra\\
          Faculty of Mathematics and Physics\\
          Sokolovsk\'a 83\\
          Charles University\\
          186 00 Praha 8\\
          Czech Republic}
 \email{tuma@karlin.mff.cuni.cz}
 \thanks{This work was completed while the first author was
 visiting the University of Caen. It was partly financed by
 the institutional grant CEZ:J13/98:113200007a and by GAUK
 no.~162/1999. The visit was financed by a \textsc{Barrande} program.}

\author[F.~Wehrung]{Friedrich Wehrung}
\address{CNRS, FRE 2271\\
D\'epartement de Math\'ematiques\\
Universit\'e de Caen\\
14032 Caen Cedex\\
France}
\email{wehrung@math.unicaen.fr}
\urladdr{http://www.math.unicaen.fr/\~{}wehrung}

\date{\today}


\begin{abstract}
Let $S$ be a distributive \jz-semilattice.
In a previous paper, the second author proved the following result:
 \begin{quote}
 \em Suppose that $S$ is a lattice.
 Let $K$ be a lattice, let $\varphi\colon\Conc K\to\nobreak S$ be a
 \jzh. Then $\varphi$ is, up to isomorphism, of the form $\Conc f$,
 for a lattice $L$ and a lattice homomorphism $f\colon K\to L$.
 \end{quote}
In the statement above, $\Conc K$ denotes as usual the \jz-semilattice
of all finitely generated congruences of $K$.

We prove here that this statement characterizes $S$
being a lattice.
\end{abstract}

\maketitle

\section*{Introduction}

The Congruence Lattice Problem (\clp\ in short) asks whether for any
distributive \jz-semilattice $S$, there exists a lattice $L$
such that $\Conc L\cong S$. While this problem is still unsolved,
many related problems have been solved. Among these, we mention the
following, due to G. Gr\"atzer and E.T. Schmidt, see
\cite{GLS1,GLS2}, and also~\cite{GrScC} for a survey about this and
related problems.

\begin{theorem}\label{T:GrSc}
Let $S$ be a finite distributive \jz-semilattice,
let $K$ be a finite lattice, let $\varphi\colon\Conc K\to\nobreak S$
be a \jzh. Then there are a finite lattice $L$, a lattice homomorphism
$f\colon K\to L$, and an isomorphism $\alpha\colon\Conc L\to S$ such
that $\alpha\circ\Conc f=\varphi$.
\end{theorem}

In the statement of Theorem~\ref{T:GrSc}, $\Conc f$ denotes the map
from $\Conc K$ to $\Conc L$ that with any congruence $\alpha$ of $K$
associates the congruence of $L$ generated by all the pairs
$\vv<f(x),f(y)>$ where $\vv<x,y>\in\alpha$.

In \cite{Wehr}, the second author proves that provided that $S$ is a
lattice, all finiteness assumptions in Theorem~\ref{T:GrSc} can be
dropped, that is:

\begin{theorem}\label{T:W}
Let $S$ be a distributive lattice with zero, let $K$ be a lattice, let
\allowbreak $\varphi\colon\Conc K\to\nobreak S$ be a \jzh. Then
$\varphi$ can be ``lifted'', that is,
there are a lattice $L$, a lattice homomorphism $f\colon K\to L$, and
an isomorphism $\alpha\colon\Conc L\to S$ such that
$\alpha\circ\Conc f=\varphi$.
\end{theorem}

In the result of Theorem~\ref{T:W}, instead of lifting a distributive
\jz-semilattice $S$ (with respect to the $\Conc$ functor), we lift a
\jzh\ $\varphi\colon\Conc K\to S$. For this reason, we shall call
such a statement ``one-dimensional Congruence Lattice Problem'', in
short $1$-\clp. With this terminology, the usual \clp\ would have to be
called $0$-\clp. By replacing $K$ by a truncated
$n$-dimensional cube (diagram) of lattices, we can define the
$n$-\clp, for any positive integer $n$. It turns out that this problem
is interesting only for
$n\in\set{0,1,2}$. Indeed, it follows from \cite{TuWe1} that the
$3$-\clp\ holds only for trivial $S$---but much more is proved in
\cite{TuWe1}, while the result about $3$-\clp\ follows from a trivial
(and unpublished) example of the second author. The
$2$-\clp\ is another matter (far less trivial than $3$-\clp\ but still
far easier than $1$-\clp), which will be considered elsewhere.

Our main result (see Theorem~A) states that for a
given distributive \jz-semilattice $S$, Theorem~\ref{T:W}
\emph{characterizes} $S$ being a lattice. This solves also
a problem formulated by H.~Dobbertin in the (yet
unpublished) monograph~\cite{Dobb}, see Corollary~\ref{C:Dobb}. In
fact, our approach is inspired by Dobbertin's solution for the
particular case of his own problem where $S$ is
\emph{primely generated}, see Theorem~15 in \cite{Dobb86}. It gives,
for a distributive \jz-semilattice $S$ that is not a lattice, the
construction of a Boolean algebra $B$ of size at most $2^{|S|}$ and a
\jzh\ $\varphi\colon\Conc B\to S$ that cannot be ``lifted'' as in
Theorem~\ref{T:W}.

Even in the particular case where $S=D$, the simplest distributive
\jz-semilattice that is not a lattice, see Section~\ref{S:Al1}, it
has been an open problem, stated at the end of Section~1 in
\cite{Dobb86}, whether the size of $B$ can be reduced from
$2^{\aleph_0}$ to $\aleph_1$ (without the Continuum Hypothesis). We
solve this affirmatively in Theorem~B. This also gives us that
\emph{there are a Boolean algebra $B$ of size $\aleph_1$ and a \jzh\
$\varphi\colon\Conc B\to D$ that cannot be lifted}, see
Corollary~\ref{C:Al1}.

We use standard notation and terminology. For a partially ordered set
$\vv<P,\leq>$ and for $a\in P$, we put
 \[
 (a]=\setm{x\in P}{x\leq a}.
 \]
We denote by $\omega$ the set of all natural numbers, and by
$\omega_1$ the first uncountable ordinal.

\section{Characterization of distributive \jz-semilattices
with $1$-\clp}\label{S:1CLP}

The main lemma of this section is the following.

\begin{lemma}\label{L:NoMD}
Let $S$ be a distributive \jz-semilattice, let $\xa_0$, $\xa_1\in S$
be such that the set $Q=(\xa_0]\cap(\xa_1]$ has no
largest element.

There are a Boolean algebra $B$
and a \jzh\ $\mu\colon B\to S$ such that the following holds:
\begin{itemize}
\item[\textup{(a)}] $\mu(1)=\xa_0\vee\xa_1$;

\item[\textup{(b)}] there are no maps $\mu_0$, $\mu_1\colon B\to S$
that satisfy the following properties:
\begin{enumerate}
\item $\mu(x)=\mu_0(x)\vee\mu_1(x)$, for all $x\in B$,

\item $\mu_0$ and $\mu_1$ are order-preserving,

\item $\mu_\ell(1)\leq\xa_\ell$, for all $\ell<2$.
\end{enumerate}
\end{itemize}

\end{lemma}

\begin{proof}
Let $\kappa$ be the minimum size of a cofinal subset of $Q$, and pick
a cofinal subset $\setm{\xx_{\xi}}{\xi<\kappa}$ of $Q$. So $\kappa$ is
an infinite cardinal. We define recursively a map
$f\colon\kappa\to\kappa$ by the rule
 \begin{equation}\label{Eq:Deff}
 f(\alpha)=\min\setm{\xi<\kappa}{\xx_{\xi}\notin
 \Id\setm{\xx_{f(\beta)}}{\beta<\alpha}}
 \end{equation}
for all $\alpha<\kappa$, where $\Id X$ denotes the ideal of $S$
generated by a subset $X$ of $S$.
Let $\beta<\alpha$. Then, by \eqref{Eq:Deff},
$\xx_{f(\alpha)}\notin\Id\setm{\xx_{f(\gamma)}}{\gamma<\alpha}$, so
$f(\alpha)\ne f(\beta)$. Moreover,
$\xx_{f(\alpha)}\notin\Id\setm{\xx_{f(\gamma)}}{\gamma<\beta}$,
so $f(\beta)\leq f(\alpha)$, whence $f(\beta)<f(\alpha)$. So $f$ is
strictly increasing.

For $\alpha<\kappa$, we put $\xq_{\alpha}=\xx_{f(\alpha)}$ and
$Q_{\alpha}=\Id\setm{\xq_{\beta}}{\beta<\alpha}$. By
\eqref{Eq:Deff},
$\xq_{\alpha}\notin Q_{\alpha}$ for all $\alpha<\kappa$. Furthermore,
all the sets $Q_{\alpha}$ are ideals of $Q$ and
$Q_{\alpha}\subset Q_{\beta}$ whenever $\alpha<\beta$. Finally, for
$\alpha<\beta$, the relation $\xq_{\alpha}\in Q_{\beta}$ holds.
(Otherwise $\xx_{f(\alpha)}\notin Q_{\beta}
=\Id\setm{\xx_{f(\gamma)}}{\gamma<\beta}$, thus, by
\eqref{Eq:Deff}, $f(\beta)\leq f(\alpha)$, a contradiction since $f$
is strictly increasing.) Hence $\bigcup_{\alpha<\kappa}Q_{\alpha}=Q$.

For $\xx\in Q$, we denote by $\|\xx\|$ the least $\alpha<\kappa$ such
that $\xx\in Q_{\alpha}$. Observe that the following obvious
properties hold:
 \begin{align}
 \|\xq_{\alpha}\|&=\alpha+1,&&\text{for all }\alpha<\kappa,
 \label{Eq:Normqa}\\
 \|\xx\vee\xy\|&=\|\xx\|\vee\|\xy\|,&&\text{for all }\xx,\,\xy\in Q.
 \label{Eq:Normxveey}
 \end{align}
Now pick a partition $\kappa=\bigcup_{\alpha<\kappa}Z_{\alpha}$ of
$\kappa$ into sets $Z_{\alpha}$ such that $|Z_{\alpha}|=\kappa$ for
all $\alpha<\kappa$. Define ideals $I$, $I_0$, and $I_1$ of the
Boolean algebra $B=P(\kappa)$ as follows:
 \begin{align*}
 I&=\setm{X\subseteq\kappa}{X\text{ finite}},\\
 I_0&=\text{ideal of }B \text{ generated by }
 \setm{Z_{\alpha}}{\alpha<\kappa},\\
 I_1&=\setm{X\subseteq\kappa}
 {X\cap Z_\alpha \text{ is finite for every }\alpha<\kappa}.
 \end{align*}
It is obvious that $I=I_0\cap I_1$, and that
$\kappa\notin I_0\cup I_1$. We define a map $\mu\colon B\to S$ by the
following rule:
 \begin{equation*}
 \mu(X)=\begin{cases}
 \bigvee_{\alpha\in X}\xq_\alpha, &\text{ if }X\text{ is finite },\\
 \xa_\ell,&\text{ if }X\in I_\ell\setminus I,\text{ for }\ell<2,\\
 \xa_0\vee\xa_1,
 &\text{ if }X\notin I_0\cup I_1.
 \end{cases}
 \end{equation*}
So $\mu$ is a \jzh\ from $B$ to $S$ with
$\mu(1)=\xa_0\vee\xa_1$.

Now suppose that $\mu_0$, $\mu_1\colon B\to S$ satisfy (i)--(iii)
above. For $\alpha<\kappa$, $\mu_1(Z_\alpha)\leq\mu(Z_\alpha)=\xa_0$
(because $Z_\alpha\in I_0\setminus I$), and
$\mu_1(Z_\alpha)\leq\mu_1(\kappa)\leq\xa_1$ (by the assumption (iii)),
hence $\mu_1(Z_\alpha)\in Q$. Hence, since $Z_{\alpha}$ is a cofinal
subset of $\kappa$, there exists $\xi_\alpha\in Z_\alpha$ such that
$\alpha\vee \|\mu_1(Z_\alpha)\|\leq\xi_\alpha$. We put
$Z=\setm{\xi_\alpha}{\alpha<\kappa}$. Observe that
$Z\in I_1\setminus I$, hence $\mu(Z)=\xa_1$. So
$\mu_0(Z)\leq\mu(Z)=\xa_1$ on the one hand, and
$\mu_0(Z)\leq\mu_0(\kappa)=\xa_0$ on the other hand, thus
$\mu_0(Z)\in Q$. Put $\beta=\|\mu_0(Z)\|$. Then
 \begin{align*}
 \xi_{\beta}+1&=\|\xq_{\xi_\beta}\|&&\text{(by \eqref{Eq:Normqa})}\\
 &=\|\mu(\set{\xi_\beta})\|&&\text{(by the definition of $\mu$)}\\
 &=\|\mu_0(\set{\xi_\beta})\|\vee\|\mu_1(\set{\xi_\beta})\|
 &&\text{(by (i) and \eqref{Eq:Normxveey})}\\
 &\leq\|\mu_0(Z)\|\vee\|\mu_1(Z_\beta)\|
 &&\text{(by (ii))}\\
 &=\beta\vee\|\mu_1(Z_\beta)\|\\
 &\leq\xi_\beta,
 \end{align*}
a contradiction.
\end{proof}

In order to formulate Corollary~\ref{C:NonMeas}, we recall the
following definition, used in particular in \cite{Weurp}. It
generalizes the classical definition of a weakly
distributive homomorphism presented in \cite{Schm68}.

\begin{definition}\label{D:WD}
Let $S$ and $T$ be join-semilattices, let $a\in S$.
A join-homomorphism $\mu\colon S\to T$ is \emph{weakly distributive}
at $a$, if for all $b_0$, $b_1\in T$ such that $\mu(a)=b_0\vee b_1$,
there are $a_0$, $a_1\in S$ such that $a=a_0\vee a_1$ and
$\mu(a_\ell)\leq b_\ell$ for all $\ell<2$.
\end{definition}

\begin{corollary}\label{C:NonMeas}
Let $S$ be a \jz-semilattice that is not a lattice.
There exist a Boolean algebra $B$ and a \jzh\
$\varphi\colon\Conc B\to S$ such that there are no lattice $L$, no
lattice homomorphism $f\colon B\to L$ and no \jzh\
$\alpha\colon\Conc L\to S$ that satisfy the following properties:
\begin{enumerate}
\item $\alpha$ is weakly distributive at $\Theta_L(f(0_B),f(1_B))$.

\item $\varphi=\alpha\circ\Conc f$.
\end{enumerate}
\end{corollary}

\begin{proof}
By assumption, there exist $\xa_0$, $\xa_1\in S$ such that
$Q=(\xa_0]\cap(\xa_1]$ has no largest element.
We consider $B$, $\mu$ as in Lemma~\ref{L:NoMD}. Since the lattice
$B$ is Boolean, the rule $x\mapsto\Theta_B(0_B,x)$ defines an
isomorphism $\pi\colon B\to\Conc B$. We put $\varphi=\mu\circ\pi^{-1}$.

So suppose that $L$, $f$, and $\alpha$ are as above. Observe that
 \[
 \alpha\Theta_L(f(0_B),f(1_B))=
 \alpha\circ(\Conc f)(\Theta_B(0_B,1_B))
 =\varphi\Theta_B(0_B,1_B)=\mu(1_B)=\xa_0\vee\xa_1,
 \]
thus, since $\alpha$ is weakly distributive at
$\Theta_L(f(0_B),f(1_B))$, there are $\Psi_0$, $\Psi_1\in\Conc L$ such
that $\Psi_0\vee\Psi_1=\Theta_L(f(0_B),f(1_B))$ and
$\alpha(\Psi_\ell)\leq\xa_\ell$, for all $\ell<2$. Thus there are
a positive integer $n$ and a decomposition
 \begin{equation}\label{Eq:Decf(0)}
 f(0_B)=t_0\leq t_1\leq\cdots\leq t_{2n}=f(1_B)
 \end{equation}
in $L$ such that the relations
 \begin{align*}
 t_{2i}&\equiv t_{2i+1}\pmod{\Psi_0},\\
 t_{2i+1}&\equiv t_{2i+2}\pmod{\Psi_1}
 \end{align*}
hold for all $i<n$. For $x\in B$, we put
 \begin{align*}
 \mu_0(x)&=\bigvee_{i<n}
 \alpha\Theta_L(t_{2i}\wedge f(x),t_{2i+1}\wedge f(x)),\\
 \mu_1(x)&=\bigvee_{i<n}
 \alpha\Theta_L(t_{2i+1}\wedge f(x),t_{2i+2}\wedge f(x)).
 \end{align*}
We verify that conditions (i)--(iii) of Lemma~\ref{L:NoMD} are
satisfied, thus causing a contradiction.
\smallskip

\noindent\textit{\textbf{Condition (i).}} For $x\in B$, we get
 \begin{align*}
 \mu_0(x)\vee\mu_1(x)&=\bigvee_{i<2n}
 \alpha\Theta_L(t_i\wedge f(x),t_{i+1}\wedge f(x))\\
 &=\alpha\Theta_L(f(0_B)\wedge f(x),f(1_B)\wedge f(x))
 &&\text{(by \eqref{Eq:Decf(0)})}\\
 &=\alpha\Theta_L(f(0_B),f(x))\\
 &=\varphi(\Theta_B(0,x))\\
 &=\mu(x).
 \end{align*}

\noindent\textit{\textbf{Condition (ii).}}
For $x\leq y$ and $i<n$, the relation
 \[
 \Theta_L(t_{2i}\wedge f(x),t_{2i+1}\wedge
 f(x))\subseteq\Theta_L(t_{2i}\wedge f(y),t_{2i+1}\wedge f(y))
 \]
holds (because $f(x)\leq f(y)$), thus $\mu_0(x)\leq\mu_0(y)$. So
$\mu_0$ is order-preserving. The proof that $\mu_1$ is
order-preserving is similar.\smallskip

\noindent\textit{\textbf{Condition (iii).}} For $i<n$,
$\Theta_L(t_{2i},t_{2i+1})\subseteq\Psi_0$, thus
$\alpha\Theta_L(t_{2i},t_{2i+1})\leq\alpha(\Psi_0)\leq\xa_0$, whence
$\mu_0(1)=\bigvee_{i<n}\alpha\Theta_L(t_{2i},t_{2i+1})\leq\xa_0$.
Similarly, $\mu_1(1)\leq\xa_1$.

This contradicts, by Lemma~\ref{L:NoMD}, the existence of $L$, $f$,
and $\alpha$.
\end{proof}

\begin{all}{Theorem A}
Let $S$ be a distributive \jz-semilattice. Then the following are
equivalent: 
\begin{enumerate}
\item For any lattice $K$ and any \jzh\
$\varphi\colon\Conc K\to S$, there are a lattice $L$, a lattice
homomorphism $f\colon K\to L$, and an isomorphism
$\alpha\colon\Conc L\to S$ such that $\varphi=\alpha\circ\Conc f$.

\item $S$ is a lattice.
\end{enumerate}
\end{all}

\begin{proof}
(ii)$\Rightarrow$(i) follows from Theorem~C in \cite{Wehr}.

(i)$\Rightarrow$(ii) is a particular case of Corollary~\ref{C:NonMeas}.
\end{proof}

With the terminology mentioned in the Introduction, this proves that
\emph{$1$-\clp\ holds at $S$ if{f} $S$ is a lattice}, for any
distributive \jz-semilattice $S$.

We also mention the following immediate consequence of
Corollary~\ref{C:NonMeas}, that solves (positively) the problem,
stated by Dobbertin in~\cite{Dobb}, whether ``strongly
measurable semilattices are lattices'':

\begin{corollary}\label{C:Dobb}
Let $S$ be a distributive \jz-semilattice. Then the following are
equivalent: 
\begin{enumerate}
\item For any Boolean algebra $B$, any \jzh\ $\mu\colon B\to S$, and
any $\xa_0$, $\xa_1\in S$ such that $\mu(1_B)=\xa_0\vee\xa_1$, there
are \jzh s $\mu_0$, $\mu_1\colon B\to S$ such that
$\mu=\mu_0\vee\mu_1$ and $\mu_\ell(1_B)=\xa_\ell$, for all $\ell<2$.

\item $S$ is a lattice.
\end{enumerate}
\end{corollary}

\begin{proof}
(ii)$\Rightarrow$(i) is proved in Corollary~10 of \cite{Dobb86}, see
also \cite{Dobb}.

(i)$\Rightarrow$(ii) follows immediately from
Corollary~\ref{C:NonMeas}.
\end{proof}

\section{A counterexample of size $\aleph_1$}\label{S:Al1}

Throughout this section, we shall denote by $D$ the \jz-semilattice
defined as $D=\omega\cup\set{\xa_0,\xa_1,\infty}$, with
$\omega$ a \jz-subsemilattice of $D$, $n<\xa_\ell<\infty$ for
all $\ell<2$, and $\infty=\xa_0\vee\xa_1$, see Figure~\ref{Fi:D}.

\begin{figure}[htb]
\begin{picture}(0,110)(40,0)
\thicklines
\put(0,0){\circle{4}}
\put(0,20){\circle{4}}
\put(0,40){\circle{4}}
\put(-20,80){\circle{4}}
\put(20,80){\circle{4}}
\put(0,100){\circle{4}}

\put(5,0){\makebox(0,0)[l]{$0$}}
\put(5,20){\makebox(0,0)[l]{$1$}}
\put(5,40){\makebox(0,0)[l]{$2$}}
\put(-25,80){\makebox(0,0)[r]{$\xa_0$}}
\put(25,80){\makebox(0,0)[l]{$\xa_1$}}
\put(0,108){\makebox(0,0){$\infty$}}

\put(0,2){\line(0,1){16}}
\put(0,22){\line(0,1){16}}
\put(-18.59,81.41){\line(1,1){17.17}}
\put(18.59,81.41){\line(-1,1){17.17}}

\multiput(0,46)(0,4){4}{\makebox(0,0){.}}
\multiput(-14.59,74.59)(4,-4){4}{\makebox(0,0){.}}
\multiput(14.59,74.59)(-4,-4){4}{\makebox(0,0){.}}

\end{picture}
\caption{The semilattice $D$}
\label{Fi:D}
\end{figure}

Now we shall construct a Boolean algebra $B$.
By Cantor's Theorem, $\aleph_1\leq 2^{\aleph_0}$,
thus there exists a one-to-one map
$f\colon\omega_1\hookrightarrow\mathcal{P}(\omega)$
(where $\mathcal{P}(\omega)$ denotes the powerset of $\omega$).
We define a map
$g\colon\omega_1\times\omega_1\to\omega$ by the rule
 \[
 g(\xi,\eta)=
 \begin{cases}
 \text{least }n<\omega\text{ such that }
 f(\xi)\cap(n+1)\ne f(\eta)\cap(n+1),&\text{if }\xi\ne\eta,\\
 0,&\text{if }\xi=\eta.
 \end{cases}
 \]

\begin{lemma}\label{L:2ton}
Let $n<\omega$, let $X$ be a subset of $\omega_1$. If
$g(\xi,\eta)<n$ for all $\xi$, $\eta\in X$, then $|X|\leq 2^n$.
\end{lemma}

\begin{proof}
Let $p$ be the map from $X$ to $\mathcal{P}(n)$ defined by the rule
 \[
 p(\xi)=f(\xi)\cap n,\quad\text{for all }\xi\in X.
 \]
(We identify $n$ with $\set{0,1,\ldots,n-1}$.)
If $|X|>2^n$, then there are $\xi$, $\eta\in X$ such that
$\xi\ne\eta$ and $p(\xi)=p(\eta)$. Hence $g(\xi,\eta)\geq n$, by the
definition of $g$, a contradiction.
\end{proof}

\begin{definition}\label{D:B}
We denote by $B$ the Boolean algebra defined by generators
$u_{0,\xi}$ and $u_{1,\xi}$, for $\xi<\omega_1$, and $v_n$, for
$n<\omega$, and the following relations:
 \begin{equation}\label{Eq:RelDefB}
 u_{0,\xi}\wedge u_{1,\eta}\leq v_{g(\xi,\eta)},\quad
 \text{for all }\xi,\,\eta<\omega_1.
 \end{equation}
Furthermore, we put $w_n=\bigvee_{k\leq n}v_k$, for all $n<\omega$.
\end{definition}

\begin{lemma}\label{L:Ineq}
$u_{0,\xi}\wedge u_{1,\eta}\leq w_n$ if{f} $g(\xi,\eta)\leq n$,
for all $\xi$, $\eta<\omega_1$ and all $n<\omega$.
\end{lemma}

\begin{proof}
If $g(\xi,\eta)\leq n$, then $u_{0,\xi}\wedge u_{1,\eta}\leq w_n$ by
\eqref{Eq:RelDefB}.

Conversely, suppose that $u_{0,\xi}\wedge u_{1,\eta}\leq w_n$. We
define elements $u_{0,\xi'}^*$, $u_{1,\eta'}^*$, and $v_k^*$ of the
two-element Boolean algebra $\two$, for $\xi'$, $\eta'<\omega_1$ and
$k<\omega$, as follows:
 \begin{align}
 u_{0,\xi}^*&=u_{1,\eta}^*=1;\label{Eq:uxieta}\\
 u_{0,\xi'}^*&=0,&&\text{for all }\xi'<\omega_1\text{ such that }
 \xi'\ne\xi;\label{Eq:u0xi'}\\
 u_{1,\eta'}^*&=0,&&\text{for all }\eta'<\omega_1\text{ such that }
 \eta'\ne\eta;\label{Eq:u1eta'}\\
 v_{g(\xi,\eta)}^*&=1;\label{Eq:vgxy}\\
 v_k^*&=0,&&\text{for all }k<\omega\text{ such that }k\ne g(\xi,\eta).
 \label{Eq:vkneg}
 \end{align}
Let $\xi'$, $\eta'<\omega_1$. If $\xi'=\xi$ and
$\eta'=\eta$, then
$u_{0,\xi'}^*\wedge u_{1,\eta'}^*=1=v_{g(\xi,\eta)}^*$. Otherwise,
$u_{0,\xi'}^*\wedge u_{1,\eta'}^*=0\leq v_{g(\xi',\eta')}^*$. So the
elements $u_{0,\xi'}^*$, $u_{1,\eta'}^*$, and $v_k^*$, for $\xi'$,
$\eta'<\omega_1$ and $k<\omega$, verify the inequalities
\eqref{Eq:RelDefB}. Therefore, there exists a homomorphism of Boolean
algebras $\varphi\colon B\to\two$ such that
 \begin{align*}
 \varphi(u_{\ell,\xi'})&=u_{\ell,\xi'}^*,&&
 \text{for all }\xi'<\omega_1\text{ and }\ell<2,\\
 \varphi(v_k)&=v_k^*,&&\text{for all }k<\omega.
 \end{align*}
In particular, by assumption,
$u_{0,\xi}^*\wedge u_{1,\eta}^*\leq\bigvee_{k\leq n}v_k^*$, that is,
$\bigvee_{k\leq n}v_k^*=1$. Therefore, by \eqref{Eq:vkneg},
$g(\xi,\eta)\leq n$.
\end{proof}

\begin{all}{Theorem B}
There exist a Boolean algebra $B$ of size $\aleph_1$ and a \jzh\
$\mu\colon B\to D$ such that the following holds:
\begin{itemize}
\item[\textup{(a)}] $\mu(1_B)=\infty$;

\item[\textup{(b)}] there are no maps $\mu_0$, $\mu_1\colon B\to D$
that satisfy the following properties:
\begin{enumerate}
\item $\mu(x)=\mu_0(x)\vee\mu_1(x)$, for all $x\in B$,

\item $\mu_0$ and $\mu_1$ are order-preserving,

\item $\mu_\ell(1)\leq\xa_\ell$, for all $\ell<2$.
\end{enumerate}
\end{itemize}
\end{all}

\begin{proof}
Let $B$ be the Boolean algebra constructed in
Definition~\ref{D:B}. It is clear that $|B|=\aleph_1$.
We define ideals $I_0$, $I_1$, and $I$ of $B$, as follows:
 \begin{align*}
 I_\ell&=\text{ideal of }B\text{ generated by }
 \setm{u_{\ell,\xi}}{\xi<\omega_1}\cup\setm{v_k}{k<\omega},
 \text{ for all }\ell<2,\\
 I&=\text{ideal of }B\text{ generated by }\setm{v_k}{k<\omega}.
  \end{align*}

It follows from \eqref{Eq:RelDefB} that $I_0\cap I_1=I$. Therefore,
we can define a \jzh\ $\mu\colon B\to D$ by the rule
 \[
 \mu(x)=
 \begin{cases}
 \text{least }n<\omega\text{ such that }x\leq w_n,&\text{if }
 x\in I_0\cap I_1,\\
 \xa_\ell,&\text{if }x\in I_\ell\setminus I_{1-\ell},
 \text{ for }\ell<2,\\
 \infty,&\text{if }x\notin I_0\cup I_1,
 \end{cases}
 \]
for all $x\in B$.

Now let $\mu_0$, $\mu_1\colon B\to D$ satisfying (i)--(iii) above.
We put
 \begin{align*}
 X_n&=\setm{\xi<\omega_1}{\mu_1(u_{0,\xi})\leq n},\\
 Y_n&=\setm{\eta<\omega_1}{\mu_0(u_{1,\eta})\leq n},
 \end{align*}
for all $n<\omega$.

\setcounter{claim}{0}
\begin{claim}\label{Cl:XYn}\hfill
\begin{itemize}
\item[\textup{(a)}] The sequences $\famm{X_n}{n<\omega}$ and
$\famm{Y_n}{n<\omega}$ are increasing.

\item[\textup{(b)}]
$\omega_1=\bigcup_{n<\omega}X_n=\bigcup_{n<\omega}Y_n$.
\end{itemize}
\end{claim}

\begin{cproof}
(i) is trivial.

(ii) Let $\xi<\omega_1$. Then
 \begin{align*}
 \mu_1(u_{0,\xi})&\leq\mu(u_{0,\xi})&&(\text{by assumption (i)})\\
 &\leq\xa_0&&(\text{by the definition of }\mu),
 \end{align*}
while also $\mu_1(u_{0,\xi})\leq\xa_1$ by assumptions (ii) and (iii).
Therefore, $\mu_1(u_{0,\xi})\leq n$ for some $n<\omega$. This proves
that $\omega_1=\bigcup_{n<\omega}X_n$. The proof that
$\omega_1=\bigcup_{n<\omega}Y_n$ is similar.
\end{cproof}

Now we put $Z_n=X_n\cap Y_n$, for all $n<\omega$. It follows from
Claim~\ref{Cl:XYn} that $\omega_1=\bigcup_{n<\omega}Z_n$. In
particular, one of the $Z_n$ should be infinite (and even
uncountable). We fix such an $n$. For all $\xi$, $\eta\in Z_n$,
$\mu_1(u_{0,\xi})\leq n$ and $\mu_0(u_{1,\eta})\leq n$, thus, by
assumptions (i) and (ii), $\mu(u_{0,\xi}\wedge u_{1,\eta})\leq n$,
that is, $u_{0,\xi}\wedge u_{1,\eta}\leq w_n$. Thus, by
Lemma~\ref{L:Ineq}, $g(\xi,\eta)\leq n$. Hence, by
Lemma~\ref{L:2ton}, $Z_n$ is finite, a contradiction.
\end{proof}

\begin{corollary}\label{C:Al1}
There exist a Boolean algebra $B$ of size $\aleph_1$ and a \jzh\
$\varphi\colon\Conc B\to D$ such that there are no lattice $L$, no
lattice homomorphism $f\colon B\to L$ and no \jzh\
$\alpha\colon\Conc L\to D$ that satisfy the following properties:
\begin{enumerate}
\item $\alpha$ is weakly distributive at $\Theta_L(f(0_B),f(1_B))$.

\item $\varphi=\alpha\circ\Conc f$.
\end{enumerate}

\end{corollary}

\begin{proof}
As in the proof of Corollary~\ref{C:NonMeas}.
\end{proof}

\section{Open problems}

The main result of Theorem~A states that the possibility, for a given
distributive \jz-semilattice $S$, to lift every \jzh\ $\Conc K\to S$
for any lattice $K$ is equivalent to $S$ being a lattice. The
maps considered in the proof of this result are not one-to-one. This
leaves open the following question:

\begin{problem}
Let $S$ be a distributive \jz-semilattice. When is it possible to
lift every \emph{one-to-one} \jzh\ 
$\varphi\colon\Conc K\hookrightarrow S$, for any lattice $K$?
\end{problem}

By Theorem~C of \cite{Wehr}, the condition that $S$ be a lattice is
sufficient. Is this condition also necessary?

\begin{problem}\label{Pb:DHom}
Let $K$ be a lattice, let $S$ be a distributive \jz-semilattice, let
$\varphi\colon\Conc K\to S$ be a \emph{distributive} \jzh. Can
$\varphi$ be lifted?
\end{problem}

Recall (see \cite{Schm68}) that for \jz-semilattices $S$ and $T$,
a homomorphism $\varphi\colon S\to T$
is \emph{distributive}, if $\varphi$ is \emph{surjective} and
$\ker\varphi$ is a directed union
of the form $\bigcup_{i\in I}\ker s_i$, where $s_i$ is a closure
operator on $S$ for all~$i$. The result of Corollary~\ref{C:NonMeas}
is of no help for solving Problem~\ref{Pb:DHom}, because the
contradiction follows there from the failure of $\alpha$ to be
(weakly) distributive.

\begin{problem}\label{Pb:AllCount}
Let $K$ be a countable lattice, let $S$ be a countable distributive
\jz-semilattice. Can every \jzh\ from $\Conc K$ to $S$ be lifted?
\end{problem}

For countable $S$, not every \jzh\ from $\Conc K$ to
$S$ can be lifted as a rule, even for $K$ of size $\aleph_1$ (this
follows from Corollary~\ref{C:Al1}). However, the problem is still open
for countable $K$.

Our last problem is more oriented to axiomatic set theory. It
originates in the observation that the construction of the Boolean
algebra of the proof of Theorem~A does not rely on
the Axiom of Choice (but it has size the continuum), while the
construction of the Boolean algebra of the proof of
Theorem~B does not rely on the Continuum Hypothesis (but it relies on
the Axiom of Choice, in the form of the existence of a one-to-one map
from $\omega_1$ into $\mathcal{P}(\omega)$).

\begin{problem}
Can one prove Theorem~B by using neither the Axiom of Choice nor the
Continuum Hypothesis?
\end{problem}

\end{document}